\documentclass[12pt,draft]{amsart}
\usepackage{amsmath,amsthm,latexsym,amscd,amsbsy,amssymb}

\chardef\bslash=`\\ 

\makeatletter
\def\verbatim{\interlinepenalty\@M \@verbatim
  \leftskip\@totalleftmargin\advance\leftskip2pc
  \frenchspacing\@vobeyspaces \@xverbatim}
\makeatother
\newtheorem{thm}{Theorem}[section]
\newtheorem{cor}[thm]{Corollary}
\newtheorem{lem}[thm]{Lemma}
\newtheorem{pro}[thm]{Proposition}

\newtheorem*{F}{Theorem (\cite[Theorem 7.3.3]{book})}
\newtheorem*{G}{Theorem (\cite[Theorem 8.1.4]{book})}
\newtheorem*{K}{Theorem ({\it \v{S}\v{c}epin}, \cite[Theorem 7.2.8]{book})}
\newtheorem*{I}{Theorem ({\it \v{S}\v{c}epin}, \cite[Theorem 8.1.6]{book})}

\theoremstyle{definition}

\theoremstyle{remark}

\numberwithin{equation}{section}

\font\f=msbm10

\begin{document}


\title[On spaces Baire isomorphic to the powers of the real line]
{On spaces Baire isomorphic to the powers of the real line}
\author{Alex Chigogidze}
\address{Department of Mathematics and Statistics,
University of Saskatche\-wan,
McLean Hall, 106 Wiggins Road, Saskatoon, SK, S7N 5E6,
Canada}
\email{chigogid@math.usask.ca}
\thanks{Author was partially supported by NSERC research grant.}
\date{\today}



\begin{abstract}{We characterize $AE(0)$-spaces that are Baire isomorphic to the powers of the real line.}
\end{abstract}

\maketitle
\markboth{A.~Chigogidze}{Baire isomorphisms}

\section{Introduction}\label{S:intro}
It is a classical result \cite{kura} that any two uncountable Polish spaces are Borel isomorphic. In particular, countable powers $\text{\f R}^{\omega}$, $\text{\f I}^{\omega}$, $\text{\f N}^{\omega}$ and $\text{\f D}^{\omega}$ of the real line $\text{\f R}$, a closed interval $\text{\f I}$, the discrete space $\text{\f N}$ of natural numbers and the two-point discrete space $\text{\f D}$ are all Borel isomorphic.

Topologically the above mentioned spaces of course differ. Characterizations of $\text{\f N}^{\omega}$ (the space of irrational numbers, \cite{alur}) and $\text{\f D}^{\omega}$ (the Cantor cube, \cite{br}) have been obtained during the early stages of the development of general topology, whereas the characterizations of spaces $\text{\f R}^{\omega}$ (the separable Hilbert space, \cite{tor2}) and $\text{\f I}^{\omega}$ (the Hilbert cube, \cite{tor1}) have been obtained relatively recently by means of powerful methods of modern infinite-dimensional topology. Later these results have been extended in order to obtain topological characterizations of {\em uncountable} powers $\text{\f R}^{\tau}$, $\text{\f I}^{\tau}$, $\text{\f N}^{\tau}$ and $\text{\f D}^{\tau}$. Turned out that characterizing properties of the above spaces, which at first glance do not have anything in common, are in fact of the same nature and all of them can be described in an unified way in terms of certain universlity properties (or, equivalently, in terms of a far going generalizations of the concept of ``general position"). 
As an illustration we recall the corresponding results.

\begin{F}\label{T:R}
Let $\tau > \omega$. The following conditions are equivalent for any $AR$-space $X$ of weight $\tau$:
\begin{enumerate}
\item
$X$ is homeomorphic to $\text{\f R}^{\tau}$.
\item
For each space $Y$ of $\text{\f R}$-weight $\leq \tau$, the set of $C$-embeddings is dense in the space $C_{\tau}(Y,X)$\footnote{Precise definition of the space $C_{\tau}(Y,X)$ is given in Section \ref{S:main}.}.
\item
For each space $Y$ of $\text{\f R}$-weight $< \tau$, the set of $C$-embeddings is dense in the space $C_{\tau}(Y,X)$.
\item
For each space $Y$ of $\text{\f R}$-weight $< \tau$ the subset
\[ \left\{ f \in C_{\tau}(Y \times \text{\f N}) \colon \text{the collection}\; \{ f(Y \times \{ n\} ) \colon n \in \text{\f N} \} \text{is discrete in}\; X\right\} \]
is dense in the space $C_{\tau}(Y \times \text{\f N})$.
\end{enumerate}
\end{F}

\begin{K}\label{T:I}
Let $\tau > \omega$. The following conditions are equivalent for any compact $AR$-space $X$ of weight $\tau$:
\begin{enumerate}
\item
$X$ is homeomorphic to $\text{\f I}^{\tau}$.
\item
For each compact space $Y$ of weight $\leq \tau$, the set of embeddings is dense in the space $C_{\tau}(Y,X)$.
\item
For each compact space $Y$ of weight $< \tau$, the set of embeddings is dense in the space $C_{\tau}(Y,X)$.
\item
For each compact space $Y$ of weight $< \tau$ the subset
\[ \left\{ f \in C_{\tau}(Y \times \text{\f D}) \colon f(Y \times \{ 0\} ) \cap f(Y\times \{ 1\} ) = \emptyset \right\} \]
is dense in the space $C_{\tau}(Y \times \text{\f D} ,X)$.
\item
$X$ is homogeneous with respect to the character, i.e. $\chi (x,X) = \tau$ for each $x \in X$.
\end{enumerate}
\end{K}

\begin{G}\label{T:N}
Let $\tau > \omega$. The following conditions are equivalent for any zero-dimensional (in the sense of $\dim$) $AE(0)$-space $X$ of weight $\tau$:
\begin{enumerate}
\item
$X$ is homeomorphic to $\text{\f N}^{\tau}$.
\item
For each zero-dimensional space $Y$ of $\text{\f R}$-weight $\leq \tau$, the set of $C$-embeddings is dense in the space $C_{\tau}(Y,X)$.
\item
For each zero-dimensional space $Y$ of $\text{\f R}$-weight $< \tau$, the set of $C$-embeddings is dense in the space $C_{\tau}(Y,X)$.
\item
For each zero-dimensional space $Y$ of $\text{\f R}$-weight $< \tau$ the subset
\[ \left\{ f \in C_{\tau}(Y \times \text{\f N}) \colon \text{the collection}\; \{ f(Y \times \{ n\} ) \colon n \in \text{\f N} \} \text{is discrete in}\; X\right\} \]
is dense in the space $C_{\tau}(Y \times \text{\f N})$.
\end{enumerate}
\end{G}

\begin{I}\label{T:D}
Let $\tau > \omega$. The following conditions are equivalent for any zero-dimensional compact $AE(0)$-space $X$ of weight $\tau$:
\begin{enumerate}
\item
$X$ is homeomorphic to $\text{\f D}^{\tau}$.
\item
For each zero-dimensional compact space $Y$ of weight $\leq \tau$, the set of embeddings is dense in the space $C_{\tau}(Y,X)$.
\item
For each zero-dimensional compact space $Y$ of weight $< \tau$, the set of embeddings is dense in the space $C_{\tau}(Y,X)$.
\item
For each zero-dimensional compact space $Y$ of weight $< \tau$ the subset
\[ \left\{ f \in C_{\tau}(Y \times \text{\f D}) \colon f(Y \times \{ 0\} ) \cap f(Y\times \{ 1\} ) = \emptyset \right\} \]
is dense in the space $C_{\tau}(Y \times \text{\f D} ,X)$.
\item
$X$ is homogeneous with respect to the character, i.e. $\chi (x,X) = \tau$ for each $x \in X$.
\end{enumerate}
\end{I}

Conditions 1--3 in the above theorems remain equivalent in the metrizable case as well and are just reformulations of the above cited results  \cite{tor2},\cite{tor1},\cite{alur} and \cite{br} respectively  (but as stated the first two results are false in the case $\tau = \omega$). Conditions 4 are variations of ``general position" properties mentioned above. It should also be emphasized that all four results are obtained by using spectral teqchnique and the general theory of $AE(n)$-spaces and $n$-soft maps. 

In the light of this discussion one would expect that the very first result cited in this article (every uncountable Polish space is Borel isomorphic to $\text{\f R}^{\omega}$) also admits a reasonable extension to the non-metrizable case and one might ask: what are spaces Baire isomorphic to $\text{\f R}^{\tau}$ for uncountable $\tau$? \footnote{Partial results in this direction are contained in \cite{sha} where the answer to the above question has been obtained for spaces representable as the limits of well-ordered spectra with perfect projections. The class of such spaces contains no uncountable power of a noncompact Polish space, let alone $\text{\f R}^{\tau}$ itself.}

Below we present a complete solution of this problem for $AE(0)$-spaces. Interestingly enough the main characterizing property (see conditions 2 and 3 of Theorem \ref{T:main}) is of the same nature as conditions 4 in the above cited results.


\section{Auxiliary results}\label{S:au}
In this section we present statements which will be used later in the proof of Theorem \ref{T:main}. All spaces are assumed to be Tychonov. Separable and completely metrizable spaces are referred as Polish spaces. Baire sets are generated by zero-sets in the same way as Borel sets by closed ones. For metrizable spaces every Borel set is a Baire set and the converse is true always. 

We assume a familiarity with the standard spectral techniques based on the \v{S}\v{c}epin's Spectral Theorem. We use \cite{book} as the main source for references, where all necessary technical details and a variety of related results can be found.

\begin{lem}\label{L:two}
Let $p \colon X \to Y$ be an open surjection of Polish spaces. Then there exists a Borel map $q \colon Y \to X$ such that $pq = \operatorname{id}_{Y}$. If, in addition, $|p^{-1}(y)| \geq 2$ for each $y \in Y$, then there exist Borel maps $q_{1}, q_{2} \colon Y \to X$ such that $pq_{i} = \operatorname{id}_{Y}$, $i = 1,2$, $q_{1}(Y) \cap q_{2}(Y) = \emptyset$ and $q_{i}(Y)$ is a Borel subset of $X$, $i = 1,2$.
\end{lem}
\begin{proof}
Let $f \colon \widetilde{Y} \to Y$ be a continuous one-to-one Borel isomorphism, where $\widetilde{Y}$ is a zero-dimensional Polish space \cite{kura}. Next consider the pullback square

\[
\begin{CD}
\widetilde{X} @>\tilde{f}>> X\\
@V{\tilde{p}}VV  @VV{p}V \\
\widetilde{Y} @>f>> Y 
\end{CD}
\]
where $\widetilde{X} = \{ (\tilde{y},x) \in \widetilde{Y} \times X \colon f(\tilde{y}) = p(x)\}$ and the maps $\tilde{p} \colon \widetilde{X} \to \widetilde{Y}$ and $\tilde{f} \colon \widetilde{X} \to X$ are the restrictions of the natural projections $\pi_{1} \colon \widetilde{Y} \times X \to \widetilde{Y}$ and $\pi_{2} \colon \widetilde{Y} \times X \to X$ onto $\widetilde{X}$ respectively.

Since the above diagram is a pullback, it follows that $\tilde{p}$ is an open surjection and $\tilde{f}$ is a contiunuous one-to-one Borel isomorphism. Since $\dim\widetilde{Y} = 0$ and since $\tilde{p}$ is open, there exists a continuous map $g_{1} \colon \widetilde{Y} \to \widetilde{X}$ such that $\tilde{p}g_{1} = \operatorname{id}_{\widetilde{Y}}$. Clearly the Borel map $q = q_{1} = \tilde{f}g_{i}f^{-1} \colon Y \to X$ satisfies the equality $qq = \operatorname{id}_{Y}$. This proves the first part of the Lemma. Obviously, $g_{1}(\widetilde{Y})$ is a closed subset of $\widetilde{X}$. Now consider a Polish space $\widetilde{X}-g_{1}(\widetilde{Y})$. Clearly the map $\tilde{p}|(\widetilde{X}-g_{1}(\widetilde{Y})) \colon \widetilde{X}-g_{1}(\widetilde{Y}) \to \widetilde{Y}$ is open and surjective (since $|\tilde{p}^{-1}(\tilde{y})| \geq 2$ for each $\tilde{y} \in \widetilde{Y}$). As above, there exists a continuous map $g_{2} \colon \widetilde{Y} \to \widetilde{X} -g_{1}(\widetilde{Y})$ such that $\tilde{p}g_{2} = \operatorname{id}_{\widetilde{Y}}$. Note that $g_{2}(\widetilde{Y})$ is also closed in $\widetilde{X}$ and $g_{1}(\widetilde{Y}) \cap g_{2}(\widetilde{Y}) = \emptyset$. Finally let $q_{i} = \tilde{f}g_{i}f^{-1} \colon Y \to X$, $i = 1,2$.
\end{proof}

\begin{lem}\label{L:mauldin}
Let $p \colon X \to Y$ be a surjective continuous map of Polish spaces. Suppose that there exists a Borel subset $M$ of $X$ such that for each $y \in Y$ the subspace $M \cap p^{-1}(y)$ is a topological copy of the Cantor cube. Then there exists a Borel isomorphism $h \colon Y \times \text{\f R}^{\omega} \to X$ such that $ph = \pi_{1}$, where $\pi_{1} \colon Y \times \text{\f R}^{\omega} \to Y$ denotes the projection onto the first coordinate.
\end{lem}
\begin{proof}
Identify $X$ with the graph of the map $p$, i.e. $X = \{ (x,y) \in X \times Y \colon p(x) = y\}$. Clearly $X$ is a closed subset of the product $X \times Y$ and the map $p$ coincides with the restriction of the projection $\pi_{1} \colon X \times Y \to Y$ onto $X$. By our assumption, there exists a Borel subset $M$ of $X$ such that $\pi_{1}^{-1}(y) \cap M$ is a copy of the Cantor set. By the main result of \cite{mauldin}, $X$ has a Borel parametrization. Since any two uncountable Polish spaces are Borel isomorphic, the latter means that there exists a Borel isomorphism $h \colon Y \times \text{\f R}^{\omega} \to X$ such that $\pi_{1}h = \pi_{1}$. Since $p = \pi_{1}|X$, it follows that $ph = \pi_{1}$ as required.
\end{proof}

\begin{lem}\label{L:pparametr}
Let ${\mathcal S} = \{ X_{n}, p_{n}^{n+1},\omega \}$ be an inverse sequence consisting of Polish spaces $X_{n}$ and open surjective projections $p_{n}^{n+1} \colon X_{n+1} \to X_{n}$. If $\left| \left( p_{n}^{n+1}\right)^{-1}(x)\right| \geq 2$ for each $n \in \omega$ and each $x \in X_{n}$, then there exists a Borel isomorphism
$h \colon X_{0} \times \text{\f R}^{\omega} \to \lim {\mathcal S}$ such that $\pi_{1} = p_{0}h$, where $p_{0} \colon \lim {\mathcal S} \to X_{0}$ is the limit projection of the spectrum ${\mathcal S}$ and $\pi_{1} \colon X_{0} \times \text{\f R}^{\omega} \to X_{0}$ is the projection onto the first coordinate.
\end{lem}
\begin{proof}
According to Lemma \ref{L:two}, for each $n \in \omega$ there exist Borel maps $q^{n}_{1}, q^{n}_{2} \colon X_{n} \to X_{n+1}$ such that $p^{n+1}_{n}q^{n}_{i} = \operatorname{id}_{X_{n}}$, $i = 1,2$, $q^{n}_{1}(X_{n}) \cap q^{n}_{2}(X_{n}) = \emptyset$ and $q^{n}_{i}(X_{n})$ is a Borel subset of $X_{n+1}$, $i = 1,2$. Let $M_{0} = X_{0}$ and suppose that for each $m \leq n$ we have already constructed a subset $M_{m} \subseteq X_{m}$ so that the following conditions are satisfied:
\begin{itemize}
\item[(a)$_{m}$]
$M_{m}$ is a Borel subset of $X_{m}$ whenever $0 \leq m \leq n$.
\item[(b)$_{m}$]
$p_{m}^{m+1}(M_{m+1}) = M_{m}$ whenever $0 \leq m \leq n-1$.
\item[(c)$_{m}$]
$\left| M_{m+1}\cap\left( p_{m}^{m+1}\right)^{-1}(x_{m})\right| = 2$ whenever $x_{m} \in M_{m}$ and $0 \leq m \leq n-1$.
\end{itemize}
We let $M_{n+1} = q_{1}^{n}(M_{n}) \cup q_{2}^{n}(M_{n})$. Since $q^{n}_{1}$ and $q^{n}_{2}$ are Borel maps, it follows that $M_{n+1}$ is a Borel subset of $X_{n+1}$. Conditions (b)$_{n+1}$ and (c)$_{n+1}$ are also satisfied by the construction. This completes the inductive step and consequently we may assume that the Borel sets $M_{n} \subseteq X_{n}$, satisfying conditions (a)$_{n}$--(c)$_{n}$, have been constructed for each $n \in \omega$. Next consider the subset $M = \cap\{ p_{n}^{-1}(M_{n}) \colon n \in \omega\}$. Conditions (a)$_{n}$, $n \in \omega$, imply that $M$ is a Borel subset of $\lim{\mathcal S}$. Conditions (b)$_{n}$ and (c)$_{n}$, $n \in \omega$, ensure that $M \cap p_{0}^{-1}(x_{0})$ is a topological copy of the Cantor cube for each $x_{0} \in X_{0}$. Finally, by Lemma \ref{L:mauldin} (with $X = \lim{\mathcal S}$, $Y = X_{0}$ and $p = p_{0}$), we conclude the existence of the required Borel isomorphism $h \colon X_{0} \times \text{\f R}^{\omega} \to \lim {\mathcal S}$.
\end{proof}

\begin{pro}\label{P:trivial}
Let ${\mathcal S} = \{ X_{n}, p_{n}^{n+1},\omega \}$ be an inverse sequence consisting of $AE(0)$-spaces $X_{n}$ and $0$-soft projections $p_{n}^{n+1} \colon X_{n+1} \to X_{n}$ which have Polish kernels. If $\left| \left( p_{n}^{n+1}\right)^{-1}(x)\right| \geq 2$ for each $n \in \omega$ and each $x \in X_{n}$, then there exists a Baire isomorphism
$h \colon X_{0} \times \text{\f R}^{\omega} \to \lim {\mathcal S}$ such that $\pi_{1} = p_{0}h$, where $p_{0} \colon \lim {\mathcal S} \to X_{0}$ is the limit projection of the spectrum ${\mathcal S}$ and $\pi_{1} \colon X_{0} \times \text{\f R}^{\omega} \to X_{0}$ is the projection onto the first coordinate.
\end{pro}
\begin{proof}
If $w(X_{0}) = \omega$, then all spaces $X_{n}$, $n \in \omega$, are Polish and the statement follows from Lemma \ref{L:pparametr}. Thus we may assume that $w(X_{0}) = \tau > \omega$. Since each short projection of the spectrum $\mathcal S$ has a Polish kernel, we conclude that $w(X_n ) = w(X_{0}) =\tau$ for each $n \in \omega$. Represent the space $X_n$, $n \in \omega$, as the limit space of a factorizing $\omega$-spectrum ${\mathcal S}_n = \{ X_{\alpha}^n ,q_{\alpha}^{\beta ,n}, A \}$, consisting of Polish spaces and $0$-soft limit projections \cite[Theorem 6.3.2]{book}. Observe, in the meantime, that the indexing sets of all these spectra coincide with $A$ (which has cardinality $\tau$). Since all short projections $p^{n+1}_n$ of the spectrum $\mathcal S$ are $0$-soft and have Polish kernels, we see, by \cite[Theorem 6.3.2(iv)]{book}, that $p^{n+1}_n$, $n \in \omega$, is the limit of some Cartesian morphism\\
$$M^{n+1}_n = \{ p^{n+1,\alpha}_n \colon X^{n+1}_{\alpha} \to X_{\alpha}^n , A_{n} \} \colon {\mathcal S}_{n+1}|A_{n} \to {\mathcal S}_{n}|A_{n} ,$$
consisting of open maps between Polish spaces, where $A_{n}$ is a cofinal and $\omega$-closed subset of the indexing set $A$.

Let $B = \cap \{ A_{n} \colon n \in \omega\}$ and note that $B$ is still a cofinal and $\omega$-closed subset of $A$ \cite[Proposition 1.1.27]{book}. In particular, $B \neq \emptyset$.

For each $\alpha \in B$ consider the inverse sequence ${\mathcal S}_{\alpha} = \{ X_{\alpha}^n , p^{n+1,\alpha}_n , \omega \}$ and let $X_{\alpha} = \lim{\mathcal S}_{\alpha}$. If $\beta \geq \alpha$, $\alpha ,\beta \in B$, then there is a Cartesian morphism\\
$$M_{\alpha}^{\beta} = \{ q_{\alpha}^{\beta ,n} \colon X_{\beta}^n \to X_{\alpha}^n ,\omega \} \colon {\mathcal S}_{\beta} \to {\mathcal S}_{\alpha} ,$$
consisting of open surjective maps $q_{\alpha}^{\beta ,n}$, $n \in \omega$. Denote by $q_{\alpha}^{\beta}$ the limit map of the morphism $M_{\alpha}^{\beta}$. Thus, the following infinite commutative diagram arises:

\begin{picture}(330,150)
\put(0,0){$X_{\alpha}$}
\put(0,130){$\lim{\mathcal S}$}
\put(115,0){$X_{\alpha}^{n+1}$}
\put(115,130){$X_{n+1}$}
\put(30,132){\vector(1,0){20}}
\put(55,132){$\dots$}
\put(70,132){\vector(1,0){40}}
\put(5,125){\vector(0,-1){115}}
\put(123,125){\vector(0,-1){115}}
\put(9,73){$q_{\alpha}$}
\put(127,73){$q_{\alpha}^{n+1}$}
\put(18,2){\vector(1,0){32}}
\put(55,2){$\dots$}
\put(70,2){\vector(1,0){40}}
\put(140,132){\vector(1,0){50}}
\put(140,2){\vector(1,0){50}}
\put(195,0){$X_{\alpha}^n$}
\put(195,130){$X_n$}
\put(210,132){\vector(1,0){40}}
\put(255,132){$\dots$}
\put(268,132){\vector(1,0){39}}
\put(210,2){\vector(1,0){40}}
\put(255,2){$\dots$}
\put(268,2){\vector(1,0){39}}
\put(200,125){\vector(0,-1){115}}
\put(205,73){$q_{\alpha}^n$}
\put(310,130){$X_0$}
\put(310,0){$X_{\alpha}^0$}
\put(315,125){\vector(0,-1){115}}
\put(320,73){$q_{\alpha}^0$}
\put(159,136){$p^{n+1}_n$}
\put(154,6){$p^{n+1,\alpha}_n$}
\end{picture}

\bigskip

\noindent Straightforward verification shows that the limit space of the spectrum ${\mathcal S} ^{\prime} = \{ X_{\alpha}, q_{\alpha}^{\beta}, A \}$ coincides with the space $X$, and that all newly arising square diagrams are also Cartesian squares. Now take an index $\alpha \in B$. Since the above mentioned diagrams are Cartesian (pullback) squares, we conclude that $\left|\left( p^{n+1,\alpha}_{n}\right)^{-1}(x)\right| \geq 2$ for each $x \in X_{\alpha}^{n}$ and $n \in \omega$. By Lemma \ref{L:pparametr}, the limit projection $p^{\alpha}_{0} \colon X_{\alpha} \to X_{\alpha}^{0}$ of the spectrum ${\mathcal S}_{\alpha}$ is Baire isomorphic to the projection $\pi_{1}^{\alpha} \colon X_{\alpha}^{0} \times \text{\f R}^{\omega} \to X_{\alpha}^{0}$, i.e. there exists a Baire isomorphism $h_{\alpha} \colon X_{\alpha}^{0} \times \text{\f R}^{\omega} \to X_{\alpha}$ such that $p^{\alpha}_{0}h_{\alpha} = \pi_{1}^{\alpha}$. The required Baire isomorphism $h \colon X_{0} \times \text{\f R}^{\omega} \to \lim{\mathcal S}$ can now be defined by letting
\[ h = q_{\alpha}^{0}\pi_{1} \triangle h_{\alpha}(q_{\alpha}^{0}\times \operatorname{id}_{\text{\f R}^{\omega}})\]
\end{proof}

\begin{pro}\label{P:selector}
Let $p \colon X \to Y$ be a $0$-soft map of $AE(0)$-spaces. Then there exists a Baire map $q \colon X \to Y$ such that $pq = \operatorname{id}_{Y}$.
\end{pro}
\begin{proof}
First let us assume that $p$ has a Polish kernel. Then, by \cite[Theorem 6.3.2(iv)]{book}, there exists pullback diagram

\[
\begin{CD}
X @>p>> Y\\
@V{f}VV  @VV{g}V \\
X_{0} @>p_{0}>> Y_{0}, 
\end{CD}
\]

\noindent where $X_{0}$ and $Y_{0}$ are Polish spaces and the map $p_{0}$ is an open surjection. By Lemma \ref{L:two}, there exists a Borel map $q_{0} \colon Y_{0} \to X_{0}$ such that $p_{0}q_{0} = \operatorname{id}_{Y_{0}}$. Then the diagonal product 
\[ q = q_{0}g \triangle \operatorname{id}_{Y} \colon Y \to X \]
is a Baire map satisfying the equality $pq = \operatorname{id}_{Y}$.

Now consider the general case. It follows from consideration in \cite[Section 6.3]{book}, that there exists a well ordered continuous inverse spectrum ${\mathcal S} = \{ X_{\alpha}, p_{\alpha}^{\alpha +1}, \tau\}$ consisting of $AE(0)$-spaces and $0$-soft short projections with Polish kernels such that $\lim{\mathcal S} = X$, $X_{0} = Y$ and the first limit projection $p_{0} \colon \lim{\mathcal S} \to X_{0}$ of the spectrum $\mathcal S$ coincides with the given map $p$. Let $q_{0} = \operatorname{id}_{X_{0}}$ and suppose that Baire for each ordinal number $\alpha < \gamma$, where $\gamma < \tau$, we have already constructed a Baire map $q_{\alpha} \colon X_{0} \to X_{\alpha}$ so that the following conditions are satisfied:
\begin{itemize}
\item
$p_{\alpha}q_{\alpha} = \operatorname{id}_{X_{0}}$ whenever $\alpha < \gamma$.
\item
$p_{\alpha}^{\beta}q_{\beta} = q_{\alpha}$ whenever $\alpha \leq \beta < \gamma$.
\item
$q_{\beta} = \triangle\{ q_{\alpha} \colon \alpha < \beta\}$ whenever $\beta < \gamma$ is a limit ordinal number.
\end{itemize}
Let us construct a Baire map $q_{\gamma} \colon X_{0} \to X_{\gamma}$. If $\gamma$ is a limit ordinal number, then let $q_{\gamma} = \triangle\{ q_{\alpha} \colon \alpha < \gamma\}$. Since the spectrum $\mathcal S$ is continuous, $q_{\gamma}$ is a well defined Baire map. If $\gamma = \alpha +1$, then according to the first part of the proof there exists a Baire map $i \colon X_{\alpha} \to X_{\alpha +1}$ (recall that the short projection $p_{\alpha}^{\alpha +1} \colon X_{\alpha +1} \to X_{\alpha}$ has a Polish kernel) such that $p_{\alpha}^{\alpha +1}i = \operatorname{id}_{X_{\alpha}}$. Then the Baire map $q_{\alpha +1} = iq_{\alpha}$ satisfies all the required properties. This completes the inductive step and hence we may assume that Baire maps $q_{\alpha} \colon X_{0} \to X_{\alpha}$ with the above properties are constructed for each $\alpha < \tau$. Let finally $q = \triangle\{ q_{\alpha} \colon \alpha < \tau\}$. Then $q \colon X_{0} \to \lim{\mathcal S}$ is a Baire map and $p_{0}q = \operatorname{id}_{X_{0}}$.
\end{proof}

\section{Main result}\label{S:main}
Let $X$ and $Y$ be arbitrary Tychonov spaces and $\tau$ be an arbitrary infinite cardinal number. Recall (see \cite[Section 6.5.1]{book}) the definition of a topology, depending on $\tau$, on the set $X^{Y}$ of all maps from $Y$ into $X$. Let $\operatorname{cov}(X)$ denote the collection of all countable functionally open covers of the space $X$. For each map $f \colon Y \to X$ the sets of the form\\
\[B(f, \{ {\mathcal U}_t \colon t \in T \} ) = \{ g \in Y^{X} \colon g\;\text{is}\; {\mathcal U}_{t}-\text{close to}\; f\;\text{for each}\; t \in T \} ,\]
where $|T| < \tau$ and ${\mathcal U}_t \in \operatorname{cov}(X)$ for each $t \in T$, are declared to be open basic neighborhoods of the point $f$ in $Y_{\tau}^{X}$. The maps, contained in the neighborhood $B(f, \{ {\mathcal U}_t \colon t \in T \} )$ are called $\{ {\mathcal U}_t \colon t \in T \}$-close to $f$.

The space $C_{\tau}(Y,X)$ mentioned in the Introduction is the subspace of the space $X^{Y}_{\tau}$ consisting of continuous maps. By ${\mathcal B}_{\tau}(Y,X)$ we denote the subspace of $X^{Y}_{\tau}$ consisting of all Baire maps of $Y$ into $X$.

Finally we say that two subsets $A$ and $B$ of a space $X$ are Baire separated if there exists a Baire subset $M$ of $X$ such that $A \subseteq M$ and $B \cap M = \emptyset$.

\begin{thm}\label{T:main}
Let $X$ be an $AE(0)$-space of weight $\tau > \omega$. Then the following conditions are equivalent:
\begin{enumerate}
\item
$X$ is Baire isomorphic to $\text{\f R}^{\tau}$.
\item
For each space $Y$ of $\text{\f R}$-weight $\leq \tau$ the set
\[ \{ f \in {\mathcal B}_{\tau}(Y \times \text{\f D}, X) \colon  f(Y \times \{ 0\} )\; \text{and}\; f(Y \times \{ 1\} )\; \text{are Baire separated}\}\]
is dense in the space ${\mathcal B}_{\tau}(Y \times \text{\f D},X)$.
\item
For each space $Y$ of $\text{\f R}$-weight $< \tau$ the set
\[ \{ f \in {\mathcal B}_{\tau}(Y \times \text{\f D}, X) \colon  f(Y \times \{ 0\} )\; \text{and}\; f(Y \times \{ 1\} )\; \text{are Baire separated}\}\]
is dense in the space ${\mathcal B}_{\tau}(Y \times \text{\f D},X)$.
\end{enumerate}
\end{thm}
\begin{proof}
(1) $\Longrightarrow$ (2). 
Let $h \colon X \to \text{\f R}^{\tau}$ be a Baire isomorphism and consider a Baire map $f \colon Y \times \text{\f D} \to X$. Take a neighbourhood $U$ of $f$ in the space ${\mathcal B}_{\tau}(Y \times \text{\f D}, X)$. By the definition of the topology ${\mathcal B}_{\tau}(Y \times \text{\f D}, X)$ there exists a collection $\{ {\mathcal U}_{t} \colon t \in T\}$ such that ${\mathcal U}_{t} \in \operatorname{cov}(X)$ for each $t \in T$, $|T| < \tau$ and 
\begin{multline*}
f \in B(f,\{ {\mathcal U}_{t}\colon t \in T\} ) =\\
 \{ g \in {\mathcal B}_{\tau}(Y \times \text{\f D},X) \colon g\;\text{is}\; {\mathcal U}_{t}-\text{close to}\; f\; \text{for each}\; t \in T\} \subseteq U .
\end{multline*}

Let $\kappa = \max\{\omega ,|T|\}$ and note that $\kappa < \tau$. Since $X$ is an $AE(0)$-space it can be represented as the limit space of a (factorizing) $\kappa$-spectrum ${\mathcal S} = \{ X_{\alpha}, p_{\alpha}^{\beta}, \exp_{\kappa}\tau\}$ consisting of $AE(0)$-spaces of weight $\kappa$ and $0$-soft limit projections (see \cite[Proposition 6.3.3]{book}). Let also ${\mathcal S}^{\prime} = \{ \left(\text{\f R}^{\omega}\right)^{\alpha}, \pi_{\alpha}^{\beta}, \exp_{\kappa}\tau \}$ be the standard $\kappa$ spectrum consisting of $\kappa$ subproducts (of $\text{\f R}^{\tau}$ and natural projections whose limit coincides with ${\f R}^{\tau}$. By the Spectral Theorem for Baire maps (see \cite[Theorem 8.8.1]{book}), we may assume without loss of generality that $h$ is the limit of a morphism
\[ h = \lim\{ h_{\alpha} \colon X_{\alpha} \to \left(\text{\f R}^{\omega}\right)^{\alpha}\}\]
consisting of Baire isomorphisms (i.e. $\pi_{\alpha}h = h_{\alpha}p_{\alpha}$ for each $\alpha \in \exp_{\kappa}\tau$).

Note that if $M$ is a Baire subset of $X$, then there exists an index $\alpha_{M} \in \exp_{\kappa}\tau$ such that $M = p_{\alpha_{M}}^{-1}(p_{\alpha_{M}}(M))$ (recall that $\mathcal S$ is a factorizing $\kappa$-spectrum) -- we say in such a case that $M$ is $\alpha_{M}$-cylindrical.  This implies that for each $t \in T$ there exists an index $\alpha_{t} \in \exp_{\kappa}\tau$ such that every element of ${\mathcal U}_{t}$ is $\alpha_{t}$-cylindrical (recall that collection ${\mathcal U}_{t}$ is countable). Finally, since $|T| \leq \kappa$ and since $\exp_{\kappa}\tau$ is a $\kappa$-complete set, it follows that there exists an index $\alpha \in \exp_{\kappa}\tau$ such that for each $t \in T$ every element of ${\mathcal U}_{t}$ is $\alpha$-cylindrical. In this situation we have 
\[ V = \{ g \in {\mathcal B}_{\tau}(Y \times \text{\f D},\text{\f R}^{\tau}) \colon p_{\alpha}g = p_{\alpha}f\} \subseteq B(f,\{ {\mathcal U}_{t}\colon t \in T\} ) \subseteq U .\]
Improtant observation here is that the set $V$ is a neighbourhood of the point $f$ in ${\mathcal B}_{\tau}(Y \times \text{\f D},X)$ (see \cite[Lemma 6.5.1]{book}).

Now consider the Baire map $hf \colon Y \times {\f D} \to {\f R}^{\tau}$ and represent it as the diagonal product
\[ hf = \pi_{\alpha}hf \triangle \pi_{\tau -\alpha}hf .\]
Since the ${\f R}$-weight of the product $Y \times {\f D}$ does not exceed $\tau$ and since $|\tau -\alpha | = \tau$, there exists a $C$-embedding $j \colon Y \times {\f D} \to \left(\text{\f R}^{\omega}\right)^{\tau -\alpha}$.

Then the diagonal product
\[ \tilde{g} = \pi_{\alpha}hf \triangle j \colon Y \times {\f D} \to {\f R}^{\tau}\]
is a Baire map such that 
\begin{itemize}
\item
$\pi_{\alpha}\tilde{g} = \pi_{\alpha}hf$.
\item
The sets $\tilde{g}(Y \times \{ 0 \} )$ and $\tilde{g}(Y \times \{ 1\} )$ are Baire separated in ${\f R}^{\tau}$.
\end{itemize}
Let now $g = h^{-1}\tilde{g} \colon Y \times \text{\f D} \to X$. It follows that the sets $g(Y \times \{ 0\} )$ and $g(Y \times \{ 1\} )$ are Baire separated in $X$. Finally observe that 
\[ p_{\alpha}g = p_{\alpha}h^{-1}\tilde{g} = h_{\alpha}^{-1} \pi_{\alpha}\tilde{g} = h_{\alpha}^{-1}\pi_{\alpha}hf = p_{\alpha}h^{-1}hf = p_{\alpha}f \]
which shows that $g \in V$. 

(2) $\Longrightarrow$ (3). Trivial.

(3) $\Longrightarrow$ (1). Embed $X$ as a $C$-embedded subspace into $\text{\f R}^{A}$, where $|A| = \tau$. By repeating the argument presented in the proof of \cite[Theorem 6.3.1]{book} we construct a collection ${\mathcal A}$ of countable subsets of the indexing set $A$ satisfying the following properties:
\begin{itemize}
\item[(a)]
The collection ${\mathcal A}$ is cofinal and $\omega$-closed in $\exp_{\omega}A$.
\item[(b)]
If $B$ is an (arbitrary) union of elements of ${\mathcal A}$, then $\pi_{B}(X)$ is a closed and $C$-embedded $AE(0)$-subspace of the product $\text{\f R}^{B}$.
\item[(c)]
If $B$ and $C$ are (arbitrary) unions of elements of $\mathcal A$ and $C \subseteq B$, then the restrictions 
\[ \pi_{C}^{B}|\pi_{B}(X) \colon \pi_{B}(X) \to \pi_{C}(X) \text{and}\; \pi_{B}|X \colon X \to \pi_{B}(X)\]
\noindent are $0$-soft.
\end{itemize}
By  using the above collection we construct a well ordered continuous inverse spectrum ${\mathcal S} = \{ X_{\alpha}, p_{\alpha}^{\alpha +1}, \tau\}$ consisting of $AE(0)$ spaces $X_{\alpha}$ and $0$-soft short projections $p_{\alpha}^{\alpha +1} \colon X_{\alpha +1} \to X_{\alpha}$ with Polish kernels so that $\lim{\mathcal S} = X$. Since $|A| = \tau$, we can write $A = \{ a_{\alpha} \colon \alpha < \tau\}$. By property (a), there exists an element $A_{0} \in {\mathcal A}$ such that $a_{0} \in A_{0}$. Let $X_{0} = \pi_{A_{0}}(X)$. Without loss of generality we may assume that $|X_{0}| > \omega$. By properties (b) and (c), $X_{0}$ is a closed subspace of $\text{\f R}^{A_{0}}$ and the map $p_{0} = \pi_{A_{0}}|X \colon X \to X_{0}$ is $0$-soft.

Suppose that for each ordinal $\alpha < \gamma$, where $\gamma < \tau$, we have already constructed a subset $A_{\alpha} \subseteq A$ as an union of less than $\tau$ elements of the collection ${\mathcal A}$ so that $\{ a_{\beta} \colon \beta < \alpha \} \subseteq A_{\alpha}$ whenever $\alpha > 0$ and all point inverses of the map
\[ p_{\alpha}^{\alpha +1} = \pi_{A_{\alpha}}^{A_{\alpha +1}}|X_{\alpha +1}\colon X_{\alpha +1} \to X_{\alpha} \]
contain at least two points. Here $X_{\alpha} = \pi_{A_{\alpha}}(X)$ and $X_{\alpha +1} = \pi_{A_{\alpha +1}}(X)$. Suppose also that the construction has been carried out in such a way that $|A_{\alpha +1}-A_{\alpha}| \leq \omega$ for each $\alpha$ with $\alpha + 1 < \gamma$. This ensures that the $0$-soft map $p_{\alpha}^{\alpha +1}$ has a Polish kernel.

If $\gamma$ is a limit ordinal number we let $A_{\gamma} = \cup\{ A_{\alpha} \colon \alpha < \gamma\}$.

If $\gamma = \alpha +1$, then consider the $0$-soft map $p_{\alpha} \colon X \to X_{\alpha}$. By Proposition \ref{P:selector}, there exists a Baire map $q \colon X_{\alpha} \to X$ such that $p_{\alpha}q = \operatorname{id}_{X_{\alpha}}$. Next consider the composition
\[ X_{\alpha} \times \text{\f D} \xrightarrow{\pi} X_{\alpha} \xrightarrow{q} X ,\]
where $\pi \colon X_{\alpha} \times \text{\f D} \to X_{\alpha}$ denotes the natural projection onto the first coordinate. According to \cite[Lemma 6.5.1]{book}, the collection of maps  $f \colon X_{\alpha} \times \text{\f D} \to X$ satisfying the equality $p_{\alpha}f = p_{\alpha}q\pi = \pi$ is a neighbourhood of the composition $q\pi$ in ${\mathcal B}_{\tau}(X_{\alpha} \times \text{\f D},X)$. Consequently, by condition (2), for at least one map $f \colon  X_{\alpha} \times \text{\f D} \to X$ with $p_{\alpha}f = \pi$ the sets $f(X_{\alpha} \times \{ 0\} )$ and $f(X_{\alpha}\times \{ 1\} )$ are Baire separated. Since $X$ is $C$-embeddwed in $\text{\f R}^{A}$, it follows that there exists a Baire subset $M$ of $\text{\f R}^{A}$ such that 
\[ f(X_{\alpha} \times \{ 0\} ) \subseteq M\;\;\text{and}\;\; f(X_{\alpha}\times \{ 1\} ) \cap M = \emptyset .\]

Choose a {\em countable} subset $B \subseteq A$ such that $M = \pi_{B}^{-1}(\pi_{B}(M))$. This obviously implies that
\[ \pi_{B}(f(X_{\alpha} \times \{ 0\} )) \cap \pi_{B}(f(X_{\alpha} \times \{ 0\} )) = \emptyset . \]
Since $p_{\alpha}f = \pi$, it follows that $B \not\subseteq A_{\alpha}$. The cofinality of the collection ${\mathcal S}$ in $\exp_{\omega}A$ allows us to find an element $\widetilde{B} \in {\mathcal A}$ such that $B \cup \{ a_{\alpha}\} \subseteq \widetilde{B}$. Clearly 
\[ \pi_{\widetilde{B}}(f(X_{\alpha} \times \{ 0\} )) \cap \pi_{\widetilde{B}}(f(X_{\alpha} \times \{ 0\} )) = \emptyset . \]
Finally let $A_{\alpha +1} = A_{\alpha} \cup \widetilde{B}$. It then follows that the $0$-soft map $p_{\alpha +1} \colon X_{\alpha +1} \to X_{\alpha}$ has a Polish kernel and all its fibers contain at least two points.

Thus the required well ordered continuous inverse spectrum ${\mathcal S} = \{ X_{\alpha}, p_{\alpha}^{\alpha +1}, \tau\}$ has been constructed so that all pont inverses of all short projections $p_{\alpha}^{\alpha +1}$ (which have Polish kernels) contain at least two points.

A straightforward transfinite induction, based on Proposition \ref{P:trivial}, shows that $X$ is Baire isomorphic to the product $X_{0} \times \left(\text{\f R}^{\omega}\right)^{\tau}$. Since $X_{0}$ is an uncountable Polish space, it is Borel isomorphic to $\text{\f R}^{\omega}$. Thus $X$ is Baire isomorphic to $\text{\f R}^{\omega} \times \left(\text{\f R}^{\omega}\right)^{\tau} \approx \text{\f R}^{\tau}$.
\end{proof}

\begin{cor}\label{C:homogeneous}
Let $X$ be an $AE(0)$-space of weight $\omega_{1}$. Then the following conditions are equivalent:
\begin{enumerate}
\item
$X$ is Baire isomorphic to $\text{\f R}^{\omega_{1}}$.
\item
For each space $Y$ of $\text{\f R}$-weight $\leq \omega_{1}$ the set
\[ \{ f \in {\mathcal B}_{\tau}(Y \times \text{\f D}, X) \colon  f(Y \times \{ 0\} )\; \text{and}\; f(Y \times \{ 1\} )\; \text{are Baire separated}\}\]
is dense in the space ${\mathcal B}_{\tau}(Y \times \text{\f D},X)$.
\item
For each Polish space $Y$  the set
\[ \{ f \in {\mathcal B}_{\tau}(Y \times \text{\f D}, X) \colon  f(Y \times \{ 0\} )\; \text{and}\; f(Y \times \{ 1\} )\; \text{are Baire separated}\}\]
is dense in the space ${\mathcal B}_{\tau}(Y \times \text{\f D},X)$.
\item
$X$ contains no $G_{\delta}$-points.
\end{enumerate}
\end{cor}

\end{document}